\documentclass{article}%
\usepackage{amsmath}
\usepackage{amsfonts}
\usepackage{amssymb}
\usepackage{graphicx}%
\providecommand{\U}[1]{\protect \rule{.1in}{.1in}}
\newtheorem{theorem}{Theorem}

\newtheorem{corollary}[theorem]{Corollary}

\newtheorem{definition}[theorem]{Definition}

\newtheorem{proposition}[theorem]{Proposition}
\newtheorem{remark}[theorem]{Remark}

\newenvironment{proof}[1][Proof]{\noindent \textbf{#1.} }{\  \rule{0.5em}{0.5em}}
\begin{document}

\title{Law of Large Numbers and Central Limit Theorem under Nonlinear Expectations}
\author{Shige PENG\\Institute of Mathematics\\Shandong University\\250100, Jinan, China\\peng@sdu.edu.cn}
\date{version February 10, 2006}
\maketitle

\section{Introduction}

The law of large numbers (LLN) and central limit theorem (CLT) are long and
widely been known as two fundamental results in probability theory.

Recently problems of model uncertainties in statistics, measures of risk and
superhedging in finance motivated us to introduce, in \cite{Peng2006a} and
\cite{Peng2006b} (see also \cite{Peng2004}, \cite{Peng2005} and references
herein), a new notion of sublinear expectation, called \textquotedblleft%
$G$-expectation\textquotedblright, and the related
\textquotedblleft$G$-normal distribution\textquotedblright \ (see
Def. \ref{Def-Gnormal}) from which we were able to define
$G$-Brownian motion as well as the corresponding stochastic
calculus. The notion of $G$-normal distribution plays the same
important rule in the theory of sublinear expectation as that of
normal distribution in the classic probability theory. It is then
natural and interesting to ask if we have the corresponding LLN and
CLT under a sublinear expectation and, in particular, if the
corresponding limit distribution of the CLT is a $G$-normal
distribution. This paper gives an affirmative answer. The proof of
our CLT is short since we borrow a deep interior estimate of fully
nonlinear PDE in \cite{WangL} which extended a profound result of
\cite{Caff1989} (see also \cite{Caff1997}) to parabolic PDEs.\ The
assumptions of our LLN and CLT can be still improved. But the
phenomenon discovered plays the same important rule in the theory of
nonlinear expectation as that of the classical LLN and CLT in
classic probability theory.

\section{Sublinear expectations\label{sec-sub}}

Let $\Omega$ be a given set and let $\mathcal{H}$ be a linear space of real
functions defined on $\Omega$ such that if $X_{1},\cdots,X_{n}\in \mathcal{H}$
then $\varphi(X_{1},\cdots,X_{n})\in \mathcal{H}$ for each $\varphi \in
C_{poly}(\mathbb{R}^{n})$ where $C_{poly}(\mathbb{R})$ denotes the space of
continuous functions with polynomial growth, i.e., there exists constants $C$
and $k\geq0$, such that $|\varphi(x)|\leq C(1+|x|^{k})$. $\mathcal{H}$ is
considered as a space of \textquotedblleft random variables\textquotedblright.

Here we use $C_{poly}(\mathbb{R}^{n})$ in our framework only for some
technique reason. In general it can be replaced by $C_{b}(\mathbb{R}^{n})$ the
space of bounded and continuous functions, by $lip_{b}(\mathbb{R}^{n})$ the
space of bounded and and Lipschitz continuous functions, or by $L^{0}%
(\mathbb{R}^{n})$ the space of Borel measurable functions.

\begin{definition}
\textbf{A sublinear expectation }$\mathbb{E}$ on $\mathcal{H}$ is a functional
$\mathcal{H}\mapsto \lbrack-\infty,\infty]$ satisfying the following
properties: for all $X,Y\in \mathcal{H}$ such that $\mathbb{E}[|X|]$,
$\mathbb{E}[|Y|]<\infty$, we have \newline \newline \textbf{(a) Monotonicity:}
if $X\geq Y$ then $\mathbb{E}[X]\geq \mathbb{E}[Y].$\newline \textbf{(b)}
\textbf{Sub-additivity (or self--dominated property):}%
\[
\mathbb{E}[X]-\mathbb{E}[Y]\leq \mathbb{E}[X-Y].
\]
\textbf{(c) Positive homogeneity: } $\mathbb{E}[\lambda X]=\lambda
\mathbb{E}[X]$,$\  \  \forall \lambda \geq0$.\newline \textbf{(d) Constant
translatability: }$\mathbb{E}[X+c]=\mathbb{E}[X]+c$.
\end{definition}

For each given $p\geq1$, we denote by $\mathcal{H}_{p}$, the collection of
$X\in \mathcal{H}$ such that $\mathbb{E}[|X|^{p}]<\infty$. It can be checked
(see \cite{Peng2006a} and \cite{Peng2006b}) that%
\[
\mathbb{E}[|X+Y|^{p}]^{1/p}\leq \mathbb{E}[|X|^{p}]^{1/p}+\mathbb{E}%
[|Y|^{p}]^{1/p}.
\]
We also have $\mathcal{H}_{q}\subseteq \mathcal{H}_{p}$ for\ $1\leq p\leq
q<\infty$ and, if $\frac{1}{p}+\frac{1}{q}=1$, then for each $X\in
\mathcal{H}_{p}$ and $Y\in \mathcal{H}_{q}$ we have $X\cdot Y\in \mathcal{H}%
_{1}$ and
\[
\mathbb{E}[|X\cdot Y|]\leq \mathbb{E}[|X|^{p}]^{1/p}\mathbb{E}[|Y|^{q}%
]^{1/q}.\
\]
It follows that $\mathcal{H}_{p}$ is a linear space and the sublinear
expectation $\mathbb{E}[\cdot]$ naturally induces a norm\ $\left \Vert
X\right \Vert _{p}:=\mathbb{E}[|X|^{p}]^{1/p}$ on $\mathcal{H}_{p}$. The
completion of $\mathcal{H}_{p}$ under this norm forms a Banach space. The
expectation $\mathbb{E}[\cdot]$ can be extended to this Banach space as well.
This extended $\mathbb{E}[\cdot]$ still satisfies the above (a)--(d). But in
this paper only the pre-Banach space $\mathcal{H}_{p}$ is involved.

\begin{proposition}
\label{Prop-X+Y}Let $X,Y\in \mathcal{H}_{1}$ be such that $\mathbb{E}%
[Y]=-\mathbb{E}[-Y]$. Then we have%
\[
\mathbb{E}[X+Y]=\mathbb{E}[X]+\mathbb{E}[Y].
\]
In particular, if $\mathbb{E}[Y]=\mathbb{E}[-Y]=0$, then $\mathbb{E}%
[X+Y]=\mathbb{E}[X]$.
\end{proposition}

\begin{proof}
It is simply because we have $\mathbb{E}[X+Y]\leq \mathbb{E}[X]+\mathbb{E}[Y]$
and
\[
\mathbb{E}[X+Y]\geq \mathbb{E}[X]-\mathbb{E}[-Y]=\mathbb{E}[X]+\mathbb{E}%
[Y]\text{.}%
\]

\end{proof}

\section{Law of Large Numbers}

\begin{theorem}
(Law of Large Numbers) Let a sequence $X_{1},X_{2},\cdots$ in $\mathcal{H}%
_{2}$ be such that
\begin{equation}
\mathbb{E}[X_{i}^{2}]=\overline{\sigma}^{2},\  \  \  \mathbb{E}[X_{i}%
X_{i+j}]=\mathbb{E}[-X_{i}X_{i+j}]=0,\  \ i,j=1,2,\cdots, \label{eq-correl}%
\end{equation}
where $\overline{\sigma}\in(0,\infty)$ is a fixed number. Then the sum
\begin{equation}
S_{n}=X_{1}+\cdots+X_{n} \label{eq-Sn}%
\end{equation}
satisfies the following law of large numbers:
\[
\  \lim_{n\rightarrow \infty}\left \Vert \frac{S_{n}}{n}\right \Vert _{2}^{2}%
=\lim_{n\rightarrow \infty}\mathbb{E}[|\frac{S_{n}}{n}|^{2}]=0.
\]
Moreover, the convergence rate is dominated by%
\[
\mathbb{E}[|\frac{S_{n}}{n}|^{2}]\leq \frac{\overline{\sigma}^{2}}{n}.\
\]

\end{theorem}

\begin{proof}
By a simple calculation, we have, using Proposition \ref{Prop-X+Y},
\begin{align*}
\mathbb{E}[|\frac{S_{n}}{n}|^{2}]  &  =\frac{1}{n^{2}}\mathbb{E}[S_{n}%
^{2}]=\frac{1}{n^{2}}\mathbb{E}[S_{n-1}^{2}+2S_{n-1}X_{n}+X_{n}^{2}]\\
&  =\frac{1}{n^{2}}\mathbb{E}[S_{n-1}^{2}+X_{n}^{2}]\leq \frac{1}{n^{2}}\{
\mathbb{E}[S_{n-1}^{2}]+\mathbb{E}[X_{n}^{2}]\} \\
&  \leq \cdots=\frac{1}{n^{2}}n\mathbb{E}[X_{1}^{2}]=\frac{\overline{\sigma
}^{2}}{n}.
\end{align*}

\end{proof}

\begin{remark}
The above condition (\ref{eq-correl}) can be easily extended to the situation
$\mathbb{E[}(X_{i}-\mu)^{2}]=\overline{\sigma}$, $\mathbb{E[}(X_{i}%
-\mu)(X_{i+j}-\mu)]=0$ and $\mathbb{E[-}(X_{i}-\mu)(X_{i+j}-\mu)]=0$, for
$i,j=1,2,\cdots$. In this case we have
\[
\lim_{n\rightarrow \infty}\mathbb{E}[|\frac{S_{n}}{n}-\mu|^{2}]=0.
\]

\end{remark}

\section{Central Limit Theorem}

We now consider a generalization of the notion of the distribution under
$\mathbb{E}$ of a random variables. To this purpose we can make a set
$\widetilde{\Omega}$ a linear space of real functions $\widetilde{\mathcal{H}%
}$ defined on $\widetilde{\Omega}$ as well as a sublinear expectation
$\widetilde{\mathbb{E}}[\cdot]$ in exact the same way as $\Omega$,
$\mathcal{H}$ and $\mathbb{E}$ defined in Section \ref{sec-sub}. We can
similarly define $\widetilde{\mathcal{H}}_{p}$ for $p\geq1$.

\begin{definition}
Two random variables, $X\in \mathcal{H}$, under $\mathbb{E}[\cdot]$ and
$Y\in \widetilde{\mathcal{H}}$ under $\widetilde{\mathbb{E}}[\cdot]$, are said
to be identically distributed if, for each $\varphi \in C_{poly}(\mathbb{R})$
such that $\varphi(X)\in \mathcal{H}_{1}$, we have $\varphi(Y)\in
\widetilde{\mathcal{H}}_{1}$ and
\[
\mathbb{E}[\varphi(X)]=\widetilde{\mathbb{E}}[\varphi(Y)].\  \
\]

\end{definition}

\begin{definition}
A random variable $X\in \mathcal{H}$ is said to be independent under
$\mathbb{E}[\cdot]$ to $Y=(Y_{1},\cdots,Y_{n})\in \mathcal{H}^{n}$ if for each
test function $\varphi \in C_{poly}(\mathbb{R}^{n+1})$ such that $\varphi
(X,Y)\in \mathcal{H}_{1}$, we have $\varphi(X,y)\in \mathcal{H}_{1}$, for each
$y\in \mathbb{R}^{n}$ and, with $\overline{\varphi}(y):=\mathbb{E}%
[\varphi(X,y)]$, we have
\[
\mathbb{E}[\varphi(X,Y)]=\mathbb{E}[\overline{\varphi}(Y)].
\]
A random variable $X\in \mathcal{H}_{2}$ is said to be weakly independent of
$Y$ if the above test functions  $\varphi$ are taken only among, instead of
$C_{poly}(\mathbb{R}^{n+1})$,%
\[
\varphi(x,y)=\psi_{0}(y)+\psi_{1}(y)x+\psi_{2}(y)x^{2},\  \  \psi_{i}\in
C_{b}(\mathbb{R}^{n}),\  \ i=1,2,3.
\]

\end{definition}

\begin{remark}
In the case of linear expectation, this notion is just the classical
independence. Note that under sublinear expectations \textquotedblleft$X$ is
independent to $Y$\textquotedblright \ does not implies automatically that
\textquotedblleft$\not Y$ is independent to $X$\textquotedblright.
\end{remark}

\begin{remark}
If we assume in the above law of large numbers that the sequence $X_{1,}%
X_{2},\cdots$ is dynamically independent and identically distributed from each
other and that $\mathbb{E}[X_{1}]=\mathbb{E}[-X_{1}]=0$, $\mathbb{E}[X_{1}%
^{2}]<\infty$. Then LLN holds.
\end{remark}

We denote by $lip_{b}(\mathbb{R})$ the collection of all uniformly Lipschitz
and bounded real functions on $\mathbb{R}$. It is a linear space.

\begin{definition}
A sequence of random variables $\left \{  \eta_{i}\right \}  _{i=1}^{\infty}$ in
$\mathcal{H}$ is said to converge in distribution under $\mathbb{E}$ if for
each $\varphi \in lip_{b}(\mathbb{R})$, $\left \{  \mathbb{E}[\varphi(\eta
_{i})]\right \}  _{i=1}^{\infty}$ converges.
\end{definition}

\begin{definition}
\label{Def-Gnormal}A random variable $\xi \in \widetilde{\mathcal{H}}$ is called
$G$-normal distributed under $\widetilde{\mathbb{E}}$, if for each $\varphi \in
lip_{b}(\mathbb{R})$, the following function defined by%
\[
u(t,x):=\mathbb{E}[\varphi(x+\sqrt{t}\xi)],\ (t,x)\in \lbrack0,\infty
)\times \mathbb{R}\
\]
is the unique (bounded and continuous) viscosity solution of the following
parabolic PDE defined on $[0,\infty)\times \mathbb{R}$:
\begin{equation}
\partial_{t}u-G(\partial_{xx}^{2}u)=0,\  \ u|_{t=0}=\varphi, \label{eq-G-heat}%
\end{equation}
where $G=G_{\underline{\sigma},\overline{\sigma}}(\alpha)$ is the
following sublinear function parameterized by $\underline{\sigma}$
and $\overline{\sigma }$ $\ $with $0\leq \underline{\sigma}\leq
\overline{\sigma}$:
\[
G(\alpha)=\frac{1}{2}(\overline{\sigma}^{2}\alpha^{+}-\underline{\sigma}%
^{2}\alpha^{-}),\  \  \alpha \in \mathbb{R},
\]
Here we denote $\alpha^{+}:=\max \{0,\alpha \}$ and $\alpha^{-}:=(-\alpha)^{+}.$
\end{definition}

\begin{remark}
A simple construction of a $G$-normal distributed random variable $\xi$ is to
take $\widetilde{\Omega}=\mathbb{R}$, $\widetilde{\mathcal{H}}=C_{poly}%
(\mathbb{R})$. The expectation $\widetilde{\mathbb{E}}$ is defined by
$\widetilde{\mathbb{E}}[\varphi]:=u^{\varphi}(1,0)$, where $u=u^{\varphi}$ is
the unique polynomial growth and continuous viscosity solution of
(\ref{eq-G-heat}) with $\varphi \in C_{poly}(\mathbb{R})=\widetilde
{\mathcal{H}}_{1}$. The $G$-normal distributed random variable is $\xi
(\omega)\equiv \omega$, $\omega \in \widetilde{\Omega}=\mathbb{R}$.
\end{remark}

Our main result is:

\begin{theorem}
(Central Limit Theorem) Let a sequence $\left \{  X_{i}\right \}  _{i=1}%
^{\infty}$ in $\mathcal{H}_{3}\mathcal{\ }$be identically distributed with
each others. We also assume that, each $X_{n+1}$ is independent (or weakly
independent) to $(X_{1},\cdots,X_{n})$ for $n=1,2,\cdots$. We assume
furthermore that
\[
\mathbb{E}[X_{1}]=\mathbb{E}[-X_{1}]=0\text{,\  \ }\mathbb{E}[X_{1}%
^{2}]=\overline{\sigma}^{2},\ -\mathbb{E}[-X_{1}^{2}]=\underline{\sigma}^{2},
\]
for some fixed numbers $0<\underline{\sigma}\leq \overline{\sigma}<\infty$.
Then, $\left \{  X_{i}\right \}  _{i=1}^{\infty}$ converges in law to the
$G$-normal distribution: for each $\varphi \in lip_{b}(\mathbb{R})$,
\begin{equation}
\lim_{n\rightarrow \infty}\mathbb{E}[\varphi(\frac{S_{n}}{\sqrt{n}%
})]=\widetilde{\mathbb{E}}[\varphi(\xi)].\label{CLT}%
\end{equation}
where $\xi$ is $G$-normal distributed under $\widetilde{\mathbb{E}}$..
\end{theorem}

\begin{proof}
For a function $\varphi \in lip_{b}(\mathbb{R})$ and a small but fixed $h>0$,
let $V$ be the unique viscosity solution of%
\begin{equation}
\partial_{t}V+G(\partial_{xx}^{2}V)=0,\ (t,x)\in \lbrack0,1+h]\times
\mathbb{R}\text{,}\  \ V|_{t=1+h}=\varphi,\label{eq-V}%
\end{equation}
We have, according to the definition of $G$-normal distribution%
\[
V(t,x)=\widetilde{\mathbb{E}}[\varphi(x+\sqrt{1+h-t}\xi)].
\]
Particularly,
\begin{equation}
V(h,0)=\widetilde{\mathbb{E}}[\varphi(\xi)],\  \ V(1+h,x)=\varphi
(x).\label{equ-0}%
\end{equation}
Since (\ref{eq-V}) is a uniformly parabolic PDE and $G$ is a convex function,
thus, by the interior regularity of $V$ (see Wang \cite{WangL}, Theorem 4.13)
we have
\[
\left \Vert V\right \Vert _{C^{1+\alpha/2,2+\alpha}([0,1]\times \mathbb{R}%
)}<\infty,\  \text{for some }\alpha \in(0,1).
\]
We set $\delta=\frac{1}{n}$ and $S_{0}=0$. Then \
\begin{align*}
&  V(1,\sqrt{\delta}S_{n})-V(0,0)=\sum_{i=0}^{n-1}\{V((i+1)\delta,\sqrt
{\delta}S_{i+1})-V(i\delta,\sqrt{\delta}S_{i})\} \\
&  =\sum_{i=0}^{n-1}\left \{  [V((i+1)\delta,\sqrt{\delta}S_{i+1}%
)-V(i\delta,\sqrt{\delta}S_{i+1})]+[V(i\delta,\sqrt{\delta}S_{i+1}%
)-V(i\delta,\sqrt{\delta}S_{i})]\right \}  \\
&  =\sum_{i=0}^{n-1}\left \{  \partial_{t}V(i\delta,\sqrt{\delta}S_{i}%
)\delta+\frac{1}{2}\partial_{xx}^{2}V(i\delta,\sqrt{\delta}S_{i})X_{i+1}%
^{2}\delta+\partial_{x}V(i\delta,\sqrt{\delta}S_{i})X_{i+1}\sqrt{\delta
}+I_{\delta}^{i}\right \}
\end{align*}
with, by Taylor's expansion,%
\begin{align*}
&  I_{\delta}^{i}=\int_{0}^{1}[\partial_{t}V((i+\beta)\delta,\sqrt{\delta
}S_{i+1})-\partial_{t}V(i\delta,\sqrt{\delta}S_{i+1})]d\beta \delta \\
&  +[\partial_{t}V(i\delta,\sqrt{\delta}S_{i+1})-\partial_{t}V(i\delta
,\sqrt{\delta}S_{i})]\delta \\
&  +\int_{0}^{1}\int_{0}^{1}[\partial_{xx}^{2}V(i\delta,\sqrt{\delta}%
S_{i}+\gamma \beta X_{i+1}\sqrt{\delta})-\partial_{xx}^{2}V(i\delta
,\sqrt{\delta}S_{i})]\beta d\beta d\gamma X_{i+1}^{2}\delta.
\end{align*}
Thus
\begin{align*}
&  \mathbb{E}[\sum_{i=0}^{n-1}\partial_{t}V(i\delta,\sqrt{\delta}S_{i}%
)\delta+\frac{1}{2}\partial_{xx}^{2}V(i\delta,\sqrt{\delta}S_{i})X_{i+1}%
^{2}\delta+\partial_{x}V(i\delta,\sqrt{\delta}S_{i})X_{i+1}\sqrt{\delta
}]-\mathbb{E}[-I_{\delta}]\\
&  \leq \mathbb{E}[V(1,\sqrt{\delta}S_{n})]-V(0,0)\\
&  \leq \mathbb{E}[\sum_{i=0}^{n-1}\partial_{t}V(i\delta,\sqrt{\delta}%
S_{i})\delta+\frac{1}{2}\partial_{xx}^{2}V(i\delta,\sqrt{\delta}S_{i}%
)X_{i+1}^{2}\delta+\partial_{x}V(i\delta,\sqrt{\delta}S_{i})X_{i+1}%
\sqrt{\delta}]+\mathbb{E}[I_{\delta}]
\end{align*}
Since $\mathbb{E}[\partial_{x}V(i\delta,\sqrt{\delta}S_{i})X_{i+1}\sqrt
{\delta}]=\mathbb{E}[-\partial_{x}V(i\delta,\sqrt{\delta}S_{i})X_{i+1}%
\sqrt{\delta}]=0$, and%
\[
\mathbb{E}[\frac{1}{2}\partial_{xx}^{2}V(i\delta,\sqrt{\delta}S_{i}%
)X_{i+1}^{2}\delta]=\mathbb{E}[G(\partial_{xx}^{2}V(i\delta,\sqrt{\delta}%
S_{i}))\delta]
\]
We have, by applying $\partial_{t}V(i\delta,\sqrt{\delta}S_{i})+\frac{1}%
{2}\partial_{xx}^{2}V(i\delta,\sqrt{\delta}S_{i})=0$,
\[
\mathbb{E}[\sum_{i=0}^{n-1}\partial_{t}V(i\delta,\sqrt{\delta}S_{i}%
)\delta+\frac{1}{2}\partial_{xx}^{2}V(i\delta,\sqrt{\delta}S_{i})X_{i+1}%
^{2}\delta+\partial_{x}V(i\delta,\sqrt{\delta}S_{i})X_{i+1}\sqrt{\delta}]=0.
\]
It then follows that%
\[
-\mathbb{E}[-\sum_{i=0}^{n-1}I_{\delta}^{i}]\leq \mathbb{E}[V(1,\sqrt{\delta
}S_{n})]-V(0,0)\leq \mathbb{E}[\sum_{i=0}^{n-1}I_{\delta}^{i}].
\]
But since both $\partial_{t}V$ and $\partial_{xx}^{2}V$ are uniformly $\alpha
$-h\"{o}lder continuous in $x$ and\ $\frac{\alpha}{2}$-h\"{o}lder continuous
in $t$ on $[0,1]\times \mathbb{R}$, we then have $|I_{\delta}^{i}|\leq
C\delta^{1+\alpha/2}[1+|X_{i+1}|+|X_{i+1}|^{2+\alpha}]$. It follows that
\[
\mathbb{E}[|I_{\delta}^{i}|]\leq C\delta^{1+\alpha/2}(1+\mathbb{E}%
[|X_{1}|^{\alpha}]+\mathbb{E}[|X_{1}|^{2+\alpha}]).
\]
Thus%
\begin{align*}
-C(\frac{1}{n})^{\alpha/2}(1+\mathbb{E}[|X_{1}|^{\alpha}+|X_{1}|^{2+\alpha}])
&  \leq \mathbb{E}[V(1,\sqrt{\delta}S_{n})]-V(0,0)\\
&  \leq C(\frac{1}{n})^{\alpha/2}(1+\mathbb{E}[|X_{1}|^{\alpha}+|X_{1}%
|^{2+\alpha}])
\end{align*}
As $n\rightarrow \infty$, we thus have
\begin{equation}
\lim_{n\rightarrow \infty}\mathbb{E}[V(1,\sqrt{\delta}S_{n}%
)]=V(0,0).\label{equ-h}%
\end{equation}
On the other hand, we have, for each $t,t^{\prime}\in \lbrack0,1+h]$ and
$x\in \mathbb{R}$,
\begin{align*}
|V(t,x)-V(t^{\prime},x)| &  =|\widetilde{\mathbb{E}}[\varphi(x+\sqrt{1+h-t}%
\xi)]-\widetilde{\mathbb{E}}[\varphi(\sqrt{1+h-t^{\prime}}\xi)]|\\
&  \leq|\widetilde{\mathbb{E}}[\varphi(x+\sqrt{1+h-t}\xi)-\varphi
(x+\sqrt{1+h-t^{\prime}}\xi)]|\\
&  \leq k_{\varphi}|(|\sqrt{1+h-t}-\sqrt{1+h-t^{\prime}}|)\widetilde
{\mathbb{E}}[|\xi|]\\
&  \leq C\sqrt{|t-t^{\prime}|},
\end{align*}
where $k_{\varphi}$ denotes the Lipschitz constant of $\varphi$. Thus
$|V(0,0)-V(0,h)|\leq C\sqrt{h}$ and, by (\ref{equ-0}),
\begin{align*}
&  |\mathbb{E}[V(1,\sqrt{\delta}S_{n})]-\mathbb{E}[\varphi(\sqrt{\delta}%
S_{n})]|\\
&  =|\mathbb{E}[V(1,\sqrt{\delta}S_{n})]-\mathbb{E}[V(1+h,\sqrt{\delta}%
S_{n})]|\leq C\sqrt{h}.
\end{align*}
It follows form (\ref{equ-h}) and (\ref{equ-0}) that
\[
\limsup_{n\rightarrow \infty}|\mathbb{E}[\varphi(\frac{S_{n}}{\sqrt{n}%
})]-\widetilde{\mathbb{E}}[\varphi(\xi)]|\leq2C\sqrt{h}.
\]
Since $h$ can be arbitrarily small we thus have%
\[
\lim_{n\rightarrow \infty}\mathbb{E}[\varphi(\frac{S_{n}}{\sqrt{n}%
})]=\widetilde{\mathbb{E}}[\varphi(\xi)].
\]

\end{proof}

\begin{corollary}
The convergence (\ref{CLT}) holds for the case where $\varphi$ is a bounded
and uniformly continuous function.
\end{corollary}

\begin{proof}
We can find a sequence $\left \{  \varphi_{k}\right \}  _{k=1}^{\infty}$ in
$lip_{b}(\mathbb{R})$ such that $\varphi_{k}\rightarrow \varphi$ uniformly on
$\mathbb{R}$. By
\begin{align*}
|\mathbb{E}[\varphi(\frac{S_{n}}{\sqrt{n}})]-\widetilde{\mathbb{E}}%
[\varphi(\xi)]|  &  \leq|\mathbb{E}[\varphi(\frac{S_{n}}{\sqrt{n}%
})]-\mathbb{E}[\varphi_{k}(\frac{S_{n}}{\sqrt{n}})]|\\
&  +|\widetilde{\mathbb{E}}[\varphi(\xi)]-\widetilde{\mathbb{E}}[\varphi
_{k}(\xi)]|+|\mathbb{E}[\varphi_{k}(\frac{S_{n}}{\sqrt{n}})]-\widetilde
{\mathbb{E}}[\varphi_{k}(\xi)]|.
\end{align*}
We can easily check that (\ref{CLT}) holds.
\end{proof}


\begin{thebibliography}{9}                                                                                                %


\bibitem {Caff1989}L.A. Caffarelli (1989) Interior estimates for fully
nonlinear equations, Ann. of Math. 130 189--213.

\bibitem {Peng2004}Peng, S. (2004) Filtration Consistent Nonlinear
Expectations and Evaluations of Contingent Claims, \emph{Acta Mathematicae
Applicatae Sinica,} English Series \textbf{20}(2), 1--24.

\bibitem {Peng2005}Peng, S. (2005) Nonlinear expectations and nonlinear Markov
chains, \emph{Chin. Ann. Math.} \textbf{26B}(2) ,159--184.

\bibitem {Peng2006a}Peng, S. (2006) $G$--Expectation, $G$--Brownian Motion and
Related Stochastic Calculus of It\^{o}'s type, preprint (pdf-file available in
arXiv:math.PR/0601035v1 3Jan 2006), to appear in \emph{Proceedings of the 2005
Abel Symposium}.

\bibitem {Peng2006b}Peng, S. (2006) Multi-Dimensional $G$-Brownian Motion and
Related Stochastic Calculus under $G$-Expectation, in arXiv:math.PR/0601699v1
28Jan 2006.

\bibitem {WangL}L. Wang, (1992) On the regularity of fully nonlinear parabolic
equations: II, Comm. Pure Appl. Math. 45, 141-178.

\bibitem {Caff1997}X. Cabre and L.A. Caffarelli, Fully nonlinear elliptic
partial di erential equations, American Math. Society (1997).
\end{thebibliography}
\end{document}